\documentclass[a4paper, 10pt, reqno]{amsart}

\usepackage{amssymb, amsmath, amsthm}
\usepackage{setspace}
\usepackage{tikz,tikz-cd}
\usepackage{mathtools} 
\usepackage{stmaryrd}

\tikzset{commutative diagrams/diagrams={baseline=-2.5pt},commutative diagrams/arrow style=tikz}
\usepackage[colorlinks]{hyperref}

\newcommand\Z{\mathbb Z}
\newcommand\C{\mathbb C}

\newcommand{\bbS}{\mathbb{S}}

\newcommand{\cC}{\mathcal{C}}
\newcommand{\cE}{\mathcal{E}}
\newcommand{\cF}{\mathcal F}

\newcommand{\cK}{\mathcal{K}}

\newcommand{\cM}{\mathcal{M}}

\newcommand{\cO}{\mathcal{O}}

\newcommand{\cT}{\mathcal T}

\newcommand{\cY}{\mathcal{Y}}

\newcommand{\mfF}{\mathfrak{F}}

\newcommand{\mfS}{\mathfrak{S}}

\newcommand{\mfW}{\mathfrak{W}}

\newcommand{\mfZ}{\mathfrak{Z}}

\newcommand\isoto{\stackrel{\sim}{\To}}

\newcommand\acts{\curvearrowright}
\newcommand\into{\hookrightarrow}

\newcommand\To{\longrightarrow}

\newcommand\Hom{\operatorname{Hom}}
\newcommand\End{\operatorname{End}}

\newcommand\Tr{\operatorname{tr}}

\renewcommand\P{\mathbb P}
\newcommand\Gr{\operatorname{Gr}}

\newcommand\rk{\operatorname{rank}}

\newcommand{\Wedge}{\mbox{\scalebox{1.2}{$\wedge$}}}

\newcommand{\GL}{\mathrm{GL}}
\newcommand\MF{\operatorname{MF}}

\newcommand{\beq}[1]{\begin{equation}\label{#1} }
\newcommand{\eeq}{\end{equation}}
\newcommand{\pgap}{\vspace{5pt}}

\theoremstyle{plain}
\newtheorem{prop}[equation]{Proposition}
\newtheorem{thm}[equation]{Theorem}
\newtheorem{lem}[equation]{Lemma}
\newtheorem{cor}[equation]{Corollary}
\newtheorem{claim}[equation]{Claim}

\theoremstyle{remark}
\newtheorem{rem}[equation]{Remark}

\theoremstyle{definition}

\makeatletter \@addtoreset{equation}{section} \makeatother

\let\oldtocsection=\tocsection
\let\oldtocsubsection=\tocsubsection
\let\oldtocsubsubsection=\tocsubsubsection
\renewcommand{\tocsection}[3]{\hspace{0em}\oldtocsection{#1}{#2}{#3}}
\renewcommand{\tocsubsection}[3]{ \hspace{1em} \oldtocsubsection{#1}{\small{#2}}{\small{#3}} }
\renewcommand{\tocsubsubsection}[3]{\hspace{2em}\oldtocsubsubsection{#1}{\small{#2}}{\small{#3}}}

\newcommand{\marginparstretch}{0.6}
\let\oldmarginpar\marginpar
\renewcommand\marginpar[1]{\-\oldmarginpar[\framebox{\setstretch{\marginparstretch}\begin{minipage}{\marginparwidth}{\raggedleft\scriptsize #1}\end{minipage}}]{\framebox{\setstretch{\marginparstretch}\begin{minipage}{\marginparwidth}{\raggedright\scriptsize #1}\end{minipage}}}}

\newcommand{\aand}{\quad\quad\mbox{and}\quad\quad}
\newcommand\ie{\emph{i.e.}~}
\newcommand{\on}{\operatorname}

\newcommand\simrightarrow{\stackrel{\textstyle\sim}{\smash{\longrightarrow}\rule{0pt}{0.3ex}}} 

\newcommand\xlra[1]{\stackrel{#1}{\longrightarrow}}

\newcommand\Fl{\mathrm{Fl}}
\newcommand\tC{\widetilde{C}}
\newcommand\Tor{\operatorname{Tor}}

\newcommand\DeltaF{\Delta}
\newcommand\cDF{\mathcal{D}}
\newcommand\cDE{\mathcal{D}'}
\newcommand\cHF{\mathcal{H}}
\newcommand\cHE{\mathcal{H}'}

\makeatletter
\newcommand{\dbloverline}[1]{\overline{\dbl@overline{#1}}}
\newcommand{\dbl@overline}[1]{\mathpalette\dbl@@overline{#1}}
\newcommand{\dbl@@overline}[2]{%
	\begingroup
	\sbox\z@{$\m@th#1\overline{#2}$}%
	\ht\z@=\dimexpr\ht\z@-2\dbl@adjust{#1}\relax
	\box\z@
	\ifx#1\scriptstyle\kern-\scriptspace\else
	\ifx#1\scriptscriptstyle\kern-\scriptspace\fi\fi
	\endgroup
}
\newcommand{\dbl@adjust}[1]{%
	\fontdimen8
	\ifx#1\displaystyle\textfont\else
	\ifx#1\textstyle\textfont\else
	\ifx#1\scriptstyle\scriptfont\else
	\scriptscriptfont\fi\fi\fi 3
}
\makeatother

\newcommand\dbl[1]{\dbloverline{#1}} 

\hypersetup{
	colorlinks=true,
	citecolor=blue,
	linkcolor=blue,
	filecolor=magenta,      
	urlcolor=cyan,
	pdftitle = {A simpler kernel for stratified Mukai flops},
	pdfauthor = {Ed Segal and Wei Tseu}
}

\begin{document}

\title{A simpler kernel for stratified Mukai flops}
\author{Ed Segal and Wei Tseu}

\begin{abstract}
We reinvestigate the problem of describing the Fourier-Mukai kernel for the derived equivalence associated to a stratified Mukai flop. We prove that if one works in the framework of matrix factorizations then the kernel has a very simple geometric construction.
\end{abstract}

\maketitle

\tableofcontents

\section{Introduction}

There is a famous conjecture of Bondal-Orlov \cite{BO} and Kawamata \cite{Kaw} that if two smooth varieties $X_+$ and $X_-$ are connected by a \emph{flop} then there should also be a derived equivalence:
$$D^b(X_+) \cong D^b(X_-)$$
Here we are using the word `flop' as a synonym for `K-equivalence', \emph{i.e.}~a birational transformation admitting a proper birational roof 
$$X_+ \stackrel{\pi_+}{\longleftarrow} \widetilde{X} \stackrel{\pi_-}{\To} X_-$$
with the property that $\pi_+^*K_{X_+} = \pi_-^* K_{X_-}$. The general case of this conjecture seems some way off, but it has been established in many classes of examples \cite{Bridgeland, Nam1, CKL, Seg, DS, Seg2, HLS, Hara, HL2, MSYZ} using a variety of different techniques. 

Some of these techniques provide an explicit Fourier-Mukai kernel 
$$\cK \in D^b(X_+\times X_-)$$
for the derived equivalence, but some do not. In the latter case, even once the equivalence has been proven, there is the follow-up problem of finding a description of the kernel. And this second problem is important; if one could formulate a plausible guess for $\cK$ in general then that would provide an entry point for proving the conjecture. 
\pgap

One guess for $\cK$ is to take the structure sheaf on the fibre product $X_+\times_{X_0} X_-$ where $X_0$ is the common contraction of $X_{\pm}$. This works for standard flops \cite{BO} and Grassmannian flops \cite{BCFMV}, where the fibre product is smooth and indeed the same as the obvious correspondence $\widetilde{X}$. It also works for Mukai flops \cite{Nam1} and Abuaf flops \cite{Hara}, where the fibre product is reducible with two components.  Unfortunately this guess fails once we get to \emph{stratified Mukai flops}
$$ X_+ = T^\vee \Gr(k,n)\; \dashleftarrow\dashrightarrow\; T^\vee \Gr(n,k) = X_- $$ 
where the fibre product has many irreducible components \cite{Nam2}.\footnote{It may be tempting to try to repair this using the derived fibre product, but this fails already for standard flops. There the derived and classical fibre products are different and it is the classical one that gives the equivalence.}
\pgap

Nevertheless, it is known that stratified Mukai flops do induce derived equivalences. This was first proven in the remarkable paper of Cautis-Kamnitzer-Licata \cite{CKL} as a consequence of a much larger structure: a categorical $\mathfrak{sl}_2$ action. A follow-up paper \cite{Cautis} provided an explicit description of the kernel $\cK$. It is a sheaf on the fibre product but not simply the structure sheaf; one must put specific line bundles on each component, away from the triple intersections, and then take a reflexive extension.\footnote{The special case where $(k,n)=(2,4)$ had been worked out earlier by Kawamata \cite{KawMukai}.}
\pgap

A very different proof of the derived equivalence for stratified Mukai flops was found by Halpern-Leistner and Sam \cite{HLS}, as a special case of a large class of examples. There are two main ingredients in their approach:
\begin{enumerate}\setlength{\itemsep}{3pt}
\item Building on work of \v{S}penko-Van den Bergh \cite{SVdB}, they use the technique of `magic windows' to prove that flops coming from quasi-symmetric GIT problems always give derived equivalences.
\item Using Kn\"orrer periodicity they replace the derived categories in the stratified Mukai flop with categories of matrix factorizations on two higher-dimensional spaces. These higher-dimensional spaces are related by a flop of the form considered in (1). 
\end{enumerate}

Recent work of the second author \cite{Tseu} shows that these two proofs produce the \emph{same} set of derived equivalences.\footnote{\emph{I.e.}~they give the same equivalence, up to shifts and the Picard group of both sides. One can make a physics argument that predicts this statement by examining the topology of the relevant stringy K\"ahler moduli space.}
\pgap

The aim of the current paper is to unpack the Halpern-Leistner--Sam equivalence and provide an explicit description of the kernel, as an object in the relevant category of matrix factorizations. We find that the answer is as simple as one could hope for; in particular it is the structure sheaf on the (reduced) fibre product of the higher-dimensional spaces.  We explain our result in Section \ref{sec.results} below. 
\pgap

Our conclusion is that the complexity of stratified Mukai flops is an artifact of Kn\"orrer periodicity. If one is happy to stay in the world of matrix factorizations then the kernel is easy to construct, and the intricacies found by Cautis \cite{Cautis} arise only when one passes back to the world of derived categories.

\subsection{Results and outline}\label{sec.results}

Let $V$ be a finite-dimensional vector space, choose $0< k < \dim V$, and let $S$ denote the tautological rank $k$ bundle on $\Gr(k, V)$. The cotangent bundle of the Grassmannian is:
$$X_+ = T^\vee \Gr(k,V) = \Hom(V/S, S)_{\Gr(k,V)} $$
It sits as a subbundle inside $\Hom(V,S)_{\Gr(k,V)}$, where it is cut out by a transverse section of the pull up of the bundle $\End(S)$. If we write $x$ and $y$ for the tautological sections of $\Hom(S, V)$ and $\Hom(V,S)$ then the section cutting out $X_+$ is:
$$y\circ x \;\in \Gamma(\End(S))$$

Now let $Q$ denote the tautological rank $k$ quotient bundle on the Grassmannian of quotients $\Gr(V,k)$. In a similar way we have a subbundle:
$$X_- = T^\vee \Gr(V,k)\;\into\;  \Hom(Q, V)_{\Gr(V,k)} $$
The \emph{Grassmannian flop} is the evident birational equivalence between the spaces $\Hom(V, S)_{\Gr(k,V)}$ and $\Hom(Q, V)_{\Gr(V,k)}$. It induces a birational equivalence between $X_+$ and $X_-$, this is the \emph{stratified Mukai flop}. 
\pgap

Now write
$$E_+ = \Hom(V,S)\oplus \End(S)_{\Gr(k,V)} $$
for the total space of the bundle $\End(S)$, pulled-up to $\Hom(V, S)_{\Gr(k,V)}$. The section $y\circ x$ becomes (since $\End(S)$ is self-dual) a function on this total space $E_+$, which we can write as
$$W = \Tr (p \circ y \circ x) \quad \in \Gamma(\cO_{E_+})$$
where $p$ is the tautological section of $\End(S)$. By Kn\"orrer periodicity \cite[...]{OrlovKP, Shipman} there is an equivalence
$$D^b(X_+) \; \cong \; \MF( E_+, W)$$
from the derived category of $X_+$ to the category of matrix factorizations of $W$. Similarly, $D^b(X_-)$ is equivalent to the category of matrix factorizations of a function $W'$ on:
$$E_- = \Hom(Q,V)\oplus \End(Q)_{\Gr(V, k)}$$
The spaces $E_+$ and $E_-$ are also birational, and the equivalence between $D^b(X_+)$ and $D^b(X_-)$ can be viewed---via Kn\"orrer periodicity---as coming from an equivalence:
$$\MF(E_+, W)\; \cong\; \MF(E_-, W')$$
In this description the Fourier-Mukai kernel for the equivalence is an object of:
$$\MF( E_+\times E_-, \; W' - W)$$

\begin{thm}\label{thm.main} Let $E^o\subset E_{\pm}$ be the open subset where $E_+$ and $E_-$ are isomorphic and write
$$\Delta\; \subset\; E^o\times E^o \;\subset\; E_+\times E_- $$
for the diagonal in this subset. Then the structure sheaf on the closure of $\Delta$
$$\cO_{\overline{\Delta}} \; \;\in \MF( E_+\times E_-, \; W' - W)$$
is the kernel for the window equivalence of \cite{HLS}.
\end{thm}
\begin{rem}\begin{enumerate}\setlength{\itemsep}{5pt}
\item The functions $W$ and $W'$ agree on $E^o$, so $W'-W$ vanishes on $\overline{\Delta}$, hence $\cO_{\overline{\Delta}}$ is indeed an object in the matrix factorization category.
\item The equivalence we are discussing restricts to the identity on $E^o$, so the kernel must \emph{a priori} be some extension of the sheaf $\cO_\Delta$. Hence our result really is the simplest we could hope for. 

\item The variety $\overline{\Delta}$ is irreducible and singular and can be described as a cone over $\Gr(k,V)\times \Gr(V,k)$. It carries a function which is the restriction of either $W$ or $W'$, and we have a correspondence of LG models:
$$\begin{tikzcd}[column sep=0em,row sep=2em]
& (\overline{\Delta},W) \arrow{dr} \arrow{dl} &   \\
(E_+, W) & \hspace{3em} & (E_-, W') 
\end{tikzcd}$$
Another way to say Theorem \ref{thm.main} is that the window equivalence is induced by this correspondence.

\item Write $E_0$ for the affine singularity underlying $E_\pm$. We prove that $\overline{\Delta}$ agrees with the reduced fibre product:
$$ \overline{\Delta} = (E_+ \times_{E_0} E_-)_{\operatorname{red}} $$
We believe that the fibre product is already reduced so in fact $\overline{\Delta}$ and $E_+\times_{E_0} E_-$ agree
 as schemes, but we have not been able to prove this.\footnote{We have verified it computationally up to $k=3$.}

\end{enumerate}
\end{rem}

We will explain the geometry of $\overline{\Delta}$ in Section \ref{sec.geometry}, and the proof of Theorem \ref{thm.main} in the remainder of Section \ref{sec.stratified}. 
\pgap

Via Kn\"orrer periodicity $\cO_{\overline{\Delta}}$ is equivalent to some object $\cK$ in $D^b(X_+\times X_-)$. By the results of \cite{Tseu}, we know that this $\cK$ must be the Fourier-Mukai kernel found by Cautis in \cite{Cautis}. In Section \ref{sec.comparing} we explain how the two constructions are related.

\subsection{AI disclosure}

ChatGPT Pro suggested Proposition \ref{prop.rationality} and Corollary \ref{cor.pushdownsplits} as a replacement for Claim \ref{claim.pushdownlocallyfree} and provided us with their proofs. Discussions with Claude and Gemini also contributed to Section \ref{sec.comparing}. The paper was written entirely by the authors. 

The first version of this paper proved only the $k=2$ case, where the key calculation (Corollary \ref{cor.weightsofC}) can be handled much more easily using a Weyman resolution. That version was written without any use of AI. 

\subsection{Acknowledgements}
Although we do not formally rely on their results we were heavily inspired by the papers \cite{BDF} and \cite{BCFMV} that discuss the problem of producing kernels for window equivalences. After completing this paper we learned that Ballard, Chidambaram and Favero have independently obtained our results in the cases $k=2$ and $3$ \cite{B}. 

We thank Calum Spicer for his generous help with Section \ref{sec.rationality}. 

W.T.\ acknowledges early support from Engineering and Physical Sciences Research Council [EP/S021590/1].

E.S. acknowledges support from OpenAI through the ChatGPT for Academic Researchers program.

\section{Magic windows and kernels}\label{sec.magicwindows}

Let $G$ be a reductive group and $U$ a representation of $G$ which is self-dual.\footnote{Or more generally quasi-symmetric.} Write $\cY$ for the Artin stack:
$$\cY = [ U/ G]$$
Any representation of $G$ defines a vector bundle on $\cY$. For a given finite set $\mfW$ of irreducible representations of $G$ we define the \emph{grade restricted subcategory}, or \emph{window}, to be the full subcategory
$$\cM_{\mfW}  \; \subset D^b(\cY) $$
generated by the corresponding finite set of vector bundles. 

In \cite{HLS} (following \cite{SVdB}) certain specific subsets $\mfW$ are defined, depending on the representation $U$, with a remarkable property.\footnote{We will not give the combinatorial details of how these sets $\mfW$ are defined but we will say explicitly what they are in the examples we care about.} The corresponding subcategories $\cM_\mfW$ are called \emph{magic windows}.

\begin{thm}\cite[Thm.~3.2]{HLS}\label{thm.magicwindow} Let $Y=U\sslash_\theta G$ be a GIT quotient of $U$ for a generic stability condition $\theta$. Then the restriction functor
$$\cM_\mfW \to D^b(Y)$$
is an equivalence.
\end{thm}
Here the GIT quotient $Y$ means, by definition, the open substack of $\cY$ consisting of $\theta$-semi-stable points. So there is a restriction functor $D^b(\cY)\to D^b(Y)$, which we are evaluating on the subcategory $\cM_\mfW$. 
\pgap

One of the key points of the above theorem is that $\mfW$ only depends on $U$ and not on $\theta$. The same window category works for all (generic) GIT quotients. Consequently if $Y_+$ and $Y_-$ are two different GIT quotients then we get a derived equivalence:
$$D^b(Y_+) \stackrel{\sim}{\longleftarrow} \cM_\mfW \stackrel{\sim}{\longrightarrow} D^b(Y_-)$$
This is called a \emph{window equivalence}. 

\begin{rem}\label{rem.GRR} Given an object $\cE \in D^b(\cY)$, its restriction $\cE\vert_0$ to the origin lies in $D^b(\mathrm{B}G)$, so is a direct sum of shifts of irreps of $G$. We will refer to these irreps as the \emph{weights} of $\cE\vert_0$. If $\cE\in \cM_\mfW$ then obviously its weights lie in $\mfW$ but the converse is generally false; for example there can exist objects supported away from the origin that do not lie in $\cM_\mfW$. However, if $\cE$ is an object which can be made equivariant for the rescaling action $\C^*\acts U$ then it has a minimal equivariant resolution, so for this class of objects we have:
$$\cE \in \cM_\mfW \;\; \iff\;\; \mbox{the weights of }\cE\vert_0 \mbox{ lie in }\mfW$$
The latter restriction is called a \emph{grade restriction rule}.
\end{rem}

Now we discuss the kernels for these window equivalences. Firstly, from Theorem \ref{thm.magicwindow} we get an embedding functor
$$D^b(Y_+) \into D^b(\cY)$$
whose image is $\cM_\mfW$. This functor has some kernel $\cK \in D^b(Y_+ \times \cY)$ with the following two properties:
\begin{enumerate}\setlength{\itemsep}{5pt}\item The restriction of $\cK$ to $Y_+\times Y_+$ is the structure sheaf $\cO_\Delta$ on the diagonal,  because the composition of the embedding functor with restriction is the identity on $D^b(Y_+)$. 

\item Consider the object:
$$\cK \vert_{Y_+\times \{0\}} \; \in D^b(Y_+\times \mathrm{B}G) $$
It has a weight decomposition, \ie it splits as a direct sum of objects in $D^b(Y_+)$ tensored with irreducible $G$-representations. Since the image of the embedding functor is $\cM_\mfW$, only irreducible representations from $\mfW$ can occur as weights of $\cK \vert_{Y_+\times \{0\}}$.
\end{enumerate}
We shall also assume:
\begin{enumerate}\item[(3)] The object $\cK$ can be made equivariant for the rescaling action $\C^*\acts Y_+\times \cY$. 
\end{enumerate}
If $G$ includes the central $\C^*\subset \GL(U)$ then this is automatic; this is not true in the examples we care about but property (3) will still be obvious.
\pgap

A consequence of Theorem \ref{thm.magicwindow} is that $\cK$ is the \emph{unique} object with these three properties; see \cite[Sect.~2.3]{HL1} and Remark \ref{rem.GRR}.\footnote{Roughly, the argument is that $D^b(Y_+\times Y_+)$ is equivalent to the subcategory $D^b(Y_+)\boxtimes \cM_\mfW$ in $D^b(Y_+\times \cY)$, and $\cK$ is the object corresponding to $\cO_\Delta$.} And once we have $\cK$ the kernel for the window equivalence is easy; it's simply the restriction of $\cK$ to the open substack $Y_+\times Y_-$. 
\pgap

So the problem is to find the object $\cK$. In principle there is an algorithm for doing this: choose an arbitrary extension $\cK'$ of $\cO_\Delta$, then there is a specific sequence of mutations which will bring it into the required window, without affecting property (1).  However in practice this is rarely feasible. A better approach is to:
\begin{itemize} \setlength{\itemsep}{3pt}\item Make the correct guess for $\cK$ the first time.
\item Prove that your guess is correct by computing the weights of $\cK \vert_{Y_+\times \{0\}}$.
\end{itemize}
This is the approach that we use in this paper. 

\begin{rem} Suppose we have a $G$-invariant function $W$ on $U$. This defines a function on $\cY$ and (by restriction) a function on any GIT quotient. We can replace derived categories with categories of matrix factorizations of $W$ in all the above, Theorem \ref{thm.magicwindow} remains true, and the story is unchanged. The function $W$ is irrelevant for the problem of finding $\cK$. 
\end{rem}

\section{Review of the Grassmannian flop}\label{sec.Gflop}

Fix vector spaces $S,V$ of dimensions $k$ and $n>k$, and consider the Artin stack:
$$ \cF \; = \; \left[  \; \Hom(S,V)\oplus \Hom(V, S) \; / \; \GL(S) \; \right] $$ 
There are two generic GIT quotients of $\cF$, associated to either positive or negative powers of the character $\det(S)$, they are:
$$F_+ = \Hom(V, S)_{\Gr(k, V)} \aand F_- = \Hom(S,V)_{\Gr(V,k)} $$
This is the Grassmannian flop. Write $F_0$ for the affine singularity that underlies them both, it is the subvariety
$$F_0 \; \subset \; \End(V)$$
of endomorphisms of rank $\leq k$. Also write $F^o \subset F_{\pm}$ for the open set where the two spaces are isomorphic. In this example three spaces coincide:
\begin{enumerate}\setlength{\itemsep}{3pt}
	\item The obvious correspondence
	$$\widetilde{F} \; = \Hom(Q, S)_{\Gr(k,V)\times \Gr(V, k)} $$
	between $F_+$ and $F_-$. 
	
	\item The fibre product $F_+\times_{F_0} F_-$.
	
	\item The closure
	$$\overline{\DeltaF} \; \subset \; F_+\times F_- $$
	where $\DeltaF\subset F^o\times F^o$ is the diagonal. In $\widetilde{F}$, the subset $\DeltaF$ is the locus where the map from $Q$ to $S$ is an isomorphism. 
	
\end{enumerate}

The structure sheaf on $\widetilde{F}$ induces a derived equivalence between $F_+$ and $F_-$ \cite{DS, BLV, BCFMV}. Let us recall the proof of this in the `magic window' framework from Section \ref{sec.magicwindows}.
\pgap

In this example a magic window can be constructed from the following set of irreducible representations of the group $\GL(S)\cong \GL_k$:
\beq{eq.window1}\mfW = \left\{ \bbS^{(a_1,\dots, a_k)}S, \;\; 0 \leq a_k \leq\dots \leq a_1 \leq n-k \, \right\}\eeq
Here $(a_1,\dots, a_k)$ is a highest weight for $\GL_k$, \ie a Young diagram, and $\bbS^{(a_1,\dots, a_k)}S$ is the corresponding Schur power of $S$, \ie the corresponding irrep. These are the irreducible representations corresponding to Kapranov's exceptional collection on $\Gr(k,n)$ \cite{Kap}. 

Now we consider the further closure:
$$ \dbl{\DeltaF} \; \subset F_+ \times \cF $$
The claim is that the structure sheaf on this locus is the kernel for the embedding functor
$$D^b(F_+) \isoto \cM_\mfW \subset D^b(\cF)$$
which implies immediately that $\cO_{\overline{\DeltaF}}$ is the kernel for the window equivalence. And as explained in the previous section, to prove this claim we just need to verify:

\begin{lem}\cite{BCFMV, Kap}\label{lem.weightsinGflop} The restriction of $\cO_{\dbl{\DeltaF}}$ to $F_+ \times \{0\}$ has weights contained in the set $\mfW$. 
\end{lem}
\begin{proof}
	We introduce coordinates on $F_+\times \cF$ as follows
	$$\begin{tikzcd}[] S_1 \ar[yshift=2pt, hook]{r}{x_1}
		&  V
		\ar[yshift=2pt]{r}{y_2}
		\ar[yshift=-2pt]{l}{y_1}
		&S_2 \ar[yshift=-2pt]{l}{x_2}\end{tikzcd}$$
	where $x_1$ is an injection. Now consider the stack
	$$ \cDF \;\; \subset \; \left[  \; \Hom(S_1,V)\oplus \Hom(V, S_2)\oplus \Hom(S_2, S_1) \; / \; \GL(S_1)\times \GL(S_2) \; \right] $$ 
	where the first map is an injection. We give $\cDF$ coordinates:
	$$\begin{tikzcd}[] S_1 \ar[yshift=2pt, hook]{r}{x_1}
		&  V
		\ar[yshift=2pt]{r}{y_2}
		&S_2 \ar[yshift=-2pt, bend left=10]{ll}{a}\end{tikzcd}$$
	There is an obvious map $\cDF \to F_+\times \cF$ by setting $y_1=ay_2$ and $x_2=x_1a$, and this map is an embedding with image  $\dbl{\DeltaF}$. Note that $\overline{\DeltaF}=\widetilde{F}$ is the locus in $\cDF$ where $y_2$ is a surjection, and $\DeltaF$ is where additionally $a$ is an isomorphism. 
	
	If we think of $\cDF$ as a vector bundle over $\Gr(k, V)\times [\Hom(V,S_2)/\GL(S_2)]$, \ie over the $x_1$ and $y_2$ coordinates, then the inclusion $\cDF\into F_+\times \cF$ is a subbundle. Let's factor this through an intermediate subbundle
	\beq{eq.Hinthemiddle} \cDF\; \into \;\cHF\; \into\; F_+\times \cF \eeq
	where $\cHF$ is the space with coordinates:
	$$\begin{tikzcd}[] S_1 \ar[yshift=2pt, hook]{r}{x_1}
		&  V
		\ar[yshift=2pt]{r}{y_2}\ar[yshift=-2pt]{l}{y_1}
		&S_2 \ar[yshift=-2pt, bend left=22]{ll}{a}\end{tikzcd} $$
	Then the quotient bundle of $\cDF$ inside $\cHF$ is $\Hom(V, S_1)$, and the quotient bundle of $\cHF$ inside $F_+\times \cF$ is 
	$$\Hom(S_2, Q_1)$$
	where $Q_1$ is the tautological quotient bundle on $\Gr(k, V)$. 
	
	Now we can compute the derived restriction of $\cO_{\cDF} = \cO_{\dbl{\DeltaF}}$ to $F_+\times \{0\}$ using a Koszul resolution. We're only interested in which $\GL(S_2)$ weights appear in the answer, and if we factor the restriction via $\cHF$ then it is only the first step that contributes to the weights. They are the $\GL(S_2)$-irreps that occur in this exterior algebra 
	$$ \Wedge^\bullet \Hom(Q_1, S_2) $$
	which is precisely the set $\mfW$. 
\end{proof}

\section{Stratified Mukai flops}\label{sec.stratified}

\subsection{The geometry of the kernel}\label{sec.geometry}

For the stratified Mukai flop, we introduce the stack
$$ \cE \; = \; \left[  \; \Hom(S,V)\oplus \Hom(V, S)\oplus \End(S) \; / \; \GL(S) \; \right] $$ 
where as before $\dim(S)=k$ and $\dim(V) = n>k$. We use coordinates $x,y$ and $p$ for the three summands. This stack contains the spaces $E_\pm$ from Section \ref{sec.results} as the open substacks
$$E_+ = \{\rk(x)=k\} \aand E_- = \{\rk(y)=k\} $$
although these are not quite the GIT quotients; see Remark \ref{rem.notGIT}. We want to understand the closure of the diagonal either in $E_+\times E_-$ or in $E_+\times \cE$.
\pgap

Recall the stack $\cDF$ from Lemma \ref{lem.weightsinGflop}, and let $\cDE$ be the total space of the vector bundle $\End(S_1)\oplus \End(S_2)$ over $\cDF$. Thus $\cDE$ has coordinates

$$\begin{tikzcd}[] S_1 \ar[yshift=2pt, hook]{r}{x_1} \ar[loop left, "p_1"]
	&  V
	\ar[yshift=2pt]{r}{y_2}
	&S_2 \ar[yshift=-2pt, bend left=10]{ll}{a}\ar[loop right, "p_2"]\end{tikzcd}$$
(where $x_1$ is an injection) and there is an embedding  $\cDE \into E_+\times \cE$. Also recall that $E^o\subset E_\pm$ denotes the common open subset, and $\Delta \subset E^o\times E^o$ the diagonal. Then $\Delta$ is a subset of $\cDE$, it is the locus where $y_2$ and $a$ have full rank and also:
$$p_2 = a^{-1}p_1 a$$ 
We can take the closure $\overline{\Delta} \subset E_+\times E_-$, or the further closure $\dbl{\Delta} \subset E_+\times \cE$, both lie inside $\cDE$.

\begin{lem}\label{lem.eqsforDeltadoublebar} The closure $\dbl{\Delta}\subset E_+\times \cE$ is the substack
	$$\dbl{\Delta} = \{\;ap_2 =  p_1 a, \;\; \chi_{p_1} = \chi_{p_2}\;\} \; \subset \cDE $$
	where $\chi_{p_i}$ is the characteristic polynomial of $p_i$.
\end{lem}
Note that the first equation alone is not enough, it cuts out a reducible subvariety with $k+1$ equidimensional components indexed by the rank of $a$; see Lemma \ref{lem.Zequidimensional}.\footnote{This is similar to the structure of the traditional fibre product $X_+\times_{X_0} X_-$ for the stratified Mukai flop, see Section \ref{sec.comparing}.}  We need the remaining equations to specify the correct irreducible component. In particular $\dbl{\Delta}$ is \emph{not} a complete intersection. 

We will prove this lemma in Section \ref{sec.weightsofC} by constructing an irreducible variety which surjects onto the locus $\{ap_2=p_1a, \chi_{p_1}=\chi_{p_2}\}$. Unfortunately we have not managed to prove that the ideal generated by the given equations is radical, so they cut out $\dbl{\Delta}$ set-theoretically but conceivably not scheme-theoretically. 
\pgap

From this we immediately get a description of $\overline{\Delta}\subset E_+\times E_-$. Take the open subset 
\beq{eq.D} D\subset \cDE\eeq
where $y_2$ is a surjection, this is a variety which is a vector bundle over $\Gr(k,V)\times \Gr(V,k)$. Then $\overline{\Delta}$ is the (reduced) cone inside this vector bundle cut out by the equations appearing in Lemma \ref{lem.eqsforDeltadoublebar}. 

It remains to compare this with the fibre product of $E_+$ and $E_-$. 

\begin{lem}\label{lem.E0} Let $E_0$ be the affine singularity underlying $\cE, E_+$ and $E_-$. Then  $E_0$ embeds into the vector space
	$$\End(V)^{\oplus k} \oplus \C^k$$
	via the map:
	$$(x,y, p) \;\mapsto\; \big(xy,\, xpy,\,\dots\, , \, xp^{k-1}y, \,  \Tr(p), \, \dots\, , \, \Tr(p^k)\big) $$
	
\end{lem}
\begin{proof}
	The claim is about the generators of the ring of $\GL(S)$-invariant functions on $\cE$. Firstly, by picking a basis for $V$ and scaling the rows of $x$ and the columns of $y$ separately we obtain a $(\C^*)^{2n}$ action on $\cE$, commuting with the $\GL(S)$ action, \ie a $2n$-dimensional multigrading on the ring of functions. By exploiting this multigrading as in \cite[\S 4] {KP} we can essentially reduce the problem to the case $n=1$.
	
	Next, using the diagonal 1-parameter subgroup in $\GL(S)$ it's clear that any invariant function must be a polynomial in the entries of $p$ and the entries of the rank 1 matrix:
	$$q = yx \quad \in \End(S)$$
	Then by \cite[\S2.5]{KP}, the Cayley-Hamilton theorem, and the fact that $\rk(q)\leq 1$, the ring of invariants (for $n=1$) is generated by:
	$$\Tr(p^s), \; \mbox{for }1\leq s \leq k  \aand \Tr(p^tq) = xp^ty,  \; \mbox{for }0\leq t < k $$
	Returning to the case of general $n$, the invariants are generated by $\Tr(p^s)$ for $s = 1,\dots, k$ (or if you prefer, by the coefficients of the characteristic polynomial of $p$) and the individual entries of the matrices $xy, xpy, \dots, xp^{k-1}y$.  
\end{proof}

\begin{cor} The fibre product $E_+\times_{E_0} E_-$ is the subscheme of $D$ cut out by the ideal $(ap_2=p_1a, \chi_{p_1}=\chi_{p_2})$, hence:
$$\overline{\Delta} = (E_+\times_{E_0} E_-)_{\operatorname{red}}$$
\end{cor}

\begin{proof} The fibre product sits inside the vector bundle $D$ \eqref{eq.D} because of the condition $x_1y_1=x_2y_2$ (as in the Grassmannian flop). From Lemma \ref{lem.E0} it is the subscheme where $p_1$ and $p_2$ have the same characteristic polynomial and also 
	$$x_1(p_1)^tay_2=x_1a(p_2)^ty_2$$
	for $t=1,\dots, k-1$. But $x_1$ is injective and $y_2$ is surjective so these equations are equivalent to $p_1a=ap_2$.
\end{proof}

As remarked in the introduction we believe that the given ideal is radical and hence the real statement is that $\overline{\Delta} = E_+\times_{E_0} E_-$. But we only have computational evidence for this stronger claim.

\subsection{The weights of the kernel}\label{sec.weights}

Now we introduce, as in Section \ref{sec.results}, the invariant function
$$W = \Tr(pyx) \quad \in \Gamma(\cO_\cE)$$
and the corresponding category of matrix factorizations $\MF(\cE, W)$. By  Kn\"orrer periodicity $\MF(\cE, W)$ is equivalent to the derived category of the zero locus:
$$\{ yx = 0 \} \quad \subset \cF $$
The two sides of the stratified Mukai flop, $X_\pm$, are the GIT quotients of this stack. This is hyperk\"ahler reduction; the equation $yx=0$ is the complex moment map for the action of $\GL(S)$ on the complex symplectic space $\Hom(S,V)\oplus \Hom(V,S)$. 

\begin{rem}\label{rem.Rcharge}
	We also include an additional $\C^*$ action on $\cE$ (the `R-charge') which acts on $(x,y,p)$ with weights $(0,0,2)$ and hence on $W$ with weight 2. This is required for the category of matrix factorizations to be $\Z$-graded and for Kn\"orrer periodicity to work as stated. But it is essentially irrelevant for our proof of Theorem \ref{thm.main} so we will avoid mentioning it as far as possible. 
\end{rem}

For the stack $\cE$ a magic window can be built from the following set of $\GL(S)$-irreps:
\begin{equation}\label{eq:mfW'}
	\mfW' =  \left\{ \bbS^{(a_1,\dots, a_k)}S, \;\; 0 \leq a_k \leq\dots \leq a_1 < n \, \right\}
\end{equation}

Note the difference between this and the set $\mfW$ for the Grassmannian flop \eqref{eq.window1}. The additional factor of $\End(S)$ relaxes the upper bound on $a_1$ (\ie on the width of the Young diagram) by $k-1$.

By Theorem \ref{thm.magicwindow} the corresponding magic window $\cM_{\mfW'} \subset D^b(\cE)$ is equivalent to the derived category of either GIT quotient. There is also a window in the matrix factorization category 
$$\cM_{\mfW'}^W \subset \MF(\cE, W) $$
defined by exactly the same grade restriction rule as $\cM_{\mfW'}$. Again the restriction functors from $\cM_{\mfW'}^W$ to matrix factorizations on either GIT quotient are equivalences \cite[Thm.~5.1]{HLS}.

\begin{rem}\label{rem.notGIT} The GIT quotients of $\cE$ are not $E_\pm$ but slightly larger open substacks. The reason is that points where $\rk x<k$ (or $\rk y<k$) can be semi-stable, provided some condition on $p$ holds. However, the critical locus of $W$ is contained entirely in $E_\pm$, so these additional semi-stable points do not affect the category of matrix factorizations. Hence the restriction functors
	$$\cM_{\mfW'}^W \to \MF(E_\pm, W)$$
	are both equivalences. This is a general fact about hyperk\"ahler quotients of this form \cite[\S 5]{HLS}.
\end{rem}
Now we consider the stack $E_+\times \cE$ with the obvious function $W_2-W_1$. The structure sheaf $\cO_{\dbl{\Delta}}$ defines an object in the category
$$\MF\!\big(E_+\times \cE, \,W_2-W_1\big)$$
since $\dbl{\Delta}$ is contained in the locus where $W_1=W_2$. Theorem \ref{thm.main} follows immediately from the following:

\begin{prop}\label{prop.weightsinSMflop} The restriction of $\cO_{\dbl{\Delta}}$ to $E_+\times \{0\}$ has weights contained in the set $\mfW'$. 
\end{prop}
\begin{proof}
	Recall from Section \ref{sec.geometry} the stack $\cDE$, which is the total space of the vector bundle $\End(S_1)\oplus \End(S_2)$ over the stack $\cDF$ that we used in Lemma \ref{lem.weightsinGflop}. Write $\cHE$ for the total space of the same vector bundle over the stack $\cHF$. Then we have inclusions:
	$$\dbl{\Delta}\;\into\; \cDE \;\into\; \cHE \;\into\;  E_+\times \cE$$
	The second and third maps here are just the pull up of \eqref{eq.Hinthemiddle}, they are inclusions of linear subbundles over the base $\Gr(k,V)\times [\Hom(V, S_2)/\GL(S_2)]$. But $\dbl{\Delta}$ is a cone in $\cDE$, cut out by the equations from Lemma \ref{lem.eqsforDeltadoublebar}. Indeed, if we define
	\beq{eq.cC}\cC = \big\{\;ap_2 =  p_1 a,\; \chi_{p_1} = \chi_{p_2}\big\} \;\; \subset \cHE \eeq
	then $\dbl{\Delta} = \cC \cap \cDE$. And since $\cDE\subset \cHE$ is the subset $\{y_1= ay_2\}$ the intersection of $\cC$ and $\cDE$ is obviously transverse.
	
	We want to compute the weights of the restriction of $\cO_{\dbl{\Delta}}$ to $E_+\times \{0\}$. We can factor the restriction via $\cHE$ just as we did in Lemma \ref{lem.weightsinGflop}. The first step is identical to the Grassmannian flop case, and contributes the weights $\mfW$. And once we are in $\cHE$ we can replace $\cO_{\dbl{\Delta}}$ with $\cO_\cC$ since they have the same weights. 
	
	Thus the problem is to compute the weights of $\cO_\cC$, \ie the additional weights that appear from the (non-transverse) equations for $\cC$.  This is a fibre-wise computation; we can ignore the $x_1, y_1, y_2$ directions and consider the cone 
	$$C \;\subset\; \big[ \Hom(S_2, S_1) \oplus \End(S_1) \oplus \End(S_2) \; / \; \GL(S_1)\times \GL(S_2)  \big]$$
	cut out by the equations \eqref{eq.cC}, then the task is to compute the $\GL(S_2)$ weights of $\cO_C$ at the origin. In Section \ref{sec.weightsofC} below we show how to complete this computation and find that the additional weights are contained in the range
	$$ 
	\bbS^{(b_1,\dots, b_k)}S_2, \quad 0 \leq b_k \leq\dots \leq b_1 < k
	$$
	(Corollary \ref{cor.weightsofC}). 	Tensoring $\mfW$ by these irreducible representations gives precisely the set $\mfW'$. 
\end{proof}

\subsection{The weights of \texorpdfstring{$C$}{C}}\label{sec.weightsofC}

Consider the affine cone:
\beq{eq.Z} Z = \{ ap_2 = p_1a\}  \;\;\subset \End(S_1)\oplus \Hom(S_2, S_1) \oplus \End(S_2) \eeq

\begin{lem}\label{lem.Zequidimensional} $Z$ is a reducible complete intersection with $k+1$ equidimensional irreducible components.
\end{lem}
\begin{proof} $Z$ is acted on by the group $\GL(S_1)\times \GL(S_2)$. If we consider the action on $\Hom(S_2, S_1)$ alone then there are finitely many orbits indexed by the possible Smith normal forms, \ie by the rank of $a$.  For a fixed $a$, the stabilizer of $a$ is the set of pairs of invertible matrices $(q_1, q_2)$ satisfying $a q_2 = q_1 a$, so this has the same dimension as the fibre of $Z$ over $a$. Therefore $Z$ is a finite union of strata which each have dimension $2k^2$.\end{proof}

As in the previous section we consider the subvariety:
$$ C = \{ap_2=p_1a,\; \chi_{p_1} = \chi_{p_2}\} \;\; \subset Z $$
We shall prove below that $C$ is irreducible, hence $C$ is the irreducible component of $Z$ given by the closure of the $\rk a=k$ locus. To complete the proof of Proposition \ref{prop.weightsinSMflop} we need to compute the weights of $\cO_C$ at the origin. 
\pgap

We will introduce a certain variety $\tC\to C$ which is generically a finite cover. To warm up we consider a similar variety that lives over $\End(S_2)$. 
\pgap

Let
$$\Fl_2\, =\, \Fl(1,2,\dots, k-1, S_2)$$
denote the full flag variety of $S_2$, and let
$$\End_\Fl (S_2) \;\subset \;\End(S_2)\times \Fl_2$$
denote the subset where the endomorphism $p_2$ preserves the flag. This is a vector bundle over $\Fl_2$ of rank $\binom{k+1}{2}$.  Note that there is a map of vector bundles $\End_\Fl (S_2) \to \C^k$ which records the eigenvalues of $p_2$ on each quotient line of the flag.

The map 
$$\psi: \End_\Fl(S_2) \To \End(S_2)$$
is generically a finite cover of degree $k!$, since if the eigenvalues of $p_2$ are distinct then a choice of preserved flag is equivalent to an ordering of the eigenvalues. But (as well as ramification) there are fibres of positive dimension, for example $\psi^{-1}(0)$ is $\Fl_2$. More complicated fibres are unions of products of flag varieties. The following statement is certainly well-known and not hard to prove, we omit the proof since we don't need it.

\begin{lem}\label{lem.easypushdownlocallyfree}
 The pushdown $\psi_*\cO$ is a $\GL(S_2)$-equivariantly trivial vector bundle of rank $k!$, and the higher derived pushdown's vanish. 
 \end{lem}

Now we turn our attention back to $C$. Write $\Fl_1$ for the flag variety of $S_1$, and consider the subset
\begin{align}\nonumber \End_\Fl(S_1) \times_{\Fl_1} \Hom_\Fl(S_2, S_1) \times_{\Fl_2}  \End_\Fl(S_2) \hspace{3cm}  \\ \label{eq.ambienttrivialbundle}
\subset \;\;\End(S_1)\times \Hom(S_2, S_1)\times \End(S_2)\times \Fl_1 \times \Fl_2 \end{align}
where $p_1, a$ and $p_2$ each preserve the flags. This is a vector bundle over $\Fl_1\times \Fl_2$ of rank  $3\binom{k+1}{2}$. Inside this vector bundle we consider a further subbundle
$$P = \left\{\begin{array}{c}p_1, a \mbox{ and } p_2 \mbox{ preserve the flags and the ordered}\\ \mbox{eigenvalues of $p_1$ and $p_2$ agree}\end{array}\right\} $$
which has rank $3\binom{k+1}{2} - k$. Finally we let:

$$\tC = P \cap  \{  ap_2 = p_1a\}$$

This space $\tC$ is not a subbundle, nor is it smooth, but it is a generalized complete intersection in $P$. The equation $ap_2-p_1a =0$ is a section of $\Hom(S_2, S_1)$, but over $P$ it is actually a section of a much smaller subbundle. Roughly speaking, at points in $P$ all  the entries of $ap_2-p_1a$ ``on or below the diagonal'' vanish automatically, and only the entries ``above the diagonal'' remain. More precisely, let $N\subset \Hom_\Fl(S_2, S_1)$ be the kernel of the map 
$$\Hom_\Fl(S_2, S_1) \to \bigoplus_{i=1}^k \Hom(L_i', L_i)$$
where $L_1,\dots, L_k$ are the tautological line bundles on $\Fl_1$ and $L_1',\dots, L_k'$ are the tautological line bundles on $\Fl_2$. 
Then over $P$, the function $ap_2-p_1a$ is a section of the bundle $N$. And if we count dimensions we see that it is a regular section. This means the space $\tC$, although singular, is amenable to computations. 

\begin{lem}\label{lem.tCirreducible} $\tC$ is irreducible. 
\end{lem}
\begin{proof} Fix a pair of flags $\mfF\in\Fl_1\times \Fl_2$ and let $\tC_\mfF$ be the fibre of $\tC$ over this point. Then we can identify $\tC_\mfF$ with the space of triples $(p_1, a, p_2)$ of upper triangular matrices such that
\beq{eq.fibreoftC} ap_2=p_1a \quad \mbox{and the diagonal entries of $p_1$ and $p_2$ agree.}\eeq
A simple count of variables and dimensions shows that every irreducible component of $\tC_\mfF$ must have dimension at least $2\binom{k+1}{2}$. 

The rest of the argument is similar to Lemma \ref{lem.Zequidimensional}. We have an action of the group of pairs of invertible upper triangular matrices, and if we consider the action on $a$ alone then there are finitely many orbits, indexed by upper triangular partial permutation matrices. If we fix an $a$ then the stabilizer of $a$ and the fibre of $\tC_\mfF$ over $a$ might not have the same dimension because of the second set of equations in \eqref{eq.fibreoftC}, but these equations are more or less redundant depending on $a$. Hence the dimension of the stratum including this $a$ lies between $ 2\binom{k+1}{2} - k$ and $2\binom{k+1}{2}$, with the upper bound attained iff $a$ is invertible. 

So $\tC_\mfF$ (and hence $\tC$) has a single irreducible component which is the closure of the locus where $a$ is invertible. 
\end{proof}

Forgetting the flags gives a map
$$\phi: \tC \to C$$
 which is generically a covering map of degree $k!$. It is easy to see that this map is surjective, hence $C$ is irreducible (as claimed above). The global version of this construction proves Lemma \ref{lem.eqsforDeltadoublebar}. 

\begin{claim}\label{claim.pushdownlocallyfree}  The (derived) pushdown $\phi_*\cO_{\tC}$ is a $\GL(S_2)$-equivariantly trivial vector bundle on $C$.
\end{claim}

If this claim holds then the weights of $\cO_{C}$ are exactly the same as the weights of $\phi_*\cO_{\tC}$, and---as we shall see below---we can calculate the latter by adjunction. 

Unfortunately, although the claim looks like an innocuous generalization of Lemma \ref{lem.easypushdownlocallyfree}, we have not been able to prove it for all $k$ (it is not hard to verify when $k=2$). Instead we shall prove a weaker statement that still suffices for our argument. The key technical point is:

\begin{prop}\label{prop.rationality}\begin{enumerate}\item[(a)] The variety $C$ has rational singularities.
\item[(b)]  It is also Gorenstein with canonical bundle $\det S_1^\vee\otimes \det S_2$.
\end{enumerate} \end{prop}

The proof of this proposition is rather involved so we leave it to the next section.  Property (a) implies:

\begin{cor}\label{cor.pushdownsplits} The unit map
$$\cO_C \to \phi_*\cO_{\tC}$$
splits in the equivariant derived category of $C$. 
\end{cor}
This may be well-known but we provide a proof for completeness.
\begin{proof}
Pick any (equivariant) resolution $C'$ of $C$ and then a resolution $\tC'$ of the main component of $\tC\times_C C'$.  So we have maps
$$\begin{tikzcd} \widetilde{C}' \ar{r}{\phi'}\ar{d}{\tilde{\rho}} & C'\ar{d}{\rho} \\ \tC \ar{r}{\phi} & C \end{tikzcd} $$
where $C'$ and $\tC'$ are smooth, $\rho, \tilde{\rho}$ are birational and $\phi, \phi'$ are generically finite covers of degree $k!$. Serre duality provides a map $\phi'_*\omega_{\phi'} \to \cO_{C'}$ which generically sums functions over the fibres. Pullback of top degree forms gives a section $\cO_{\tC'} \to \omega_{\phi'}$, and we have the unit morphism $\cO_{\tC}\to \tilde{\rho}_*\cO_{\tC'}$. Combining these three gives a map
$$\Sigma\!: \phi_*\cO_{\tC} \To \rho_* \cO_{C'}\cong \cO_C $$
since $C$ has rational singularities. Composing $\Sigma$ with the unit $\cO_C \to \phi_*\cO_{\tC}$ gives an endomorphism of $\cO_C$ which generically is multiplication by $k!$, so $\tfrac{1}{k!}\Sigma$ is the required splitting. 
\end{proof}

Hence $\cO_C$ is a direct summand of $\phi_*\cO_{\tC}$, so the weights of $\cO_C$ are a subset of the weights of $\phi_*\cO_{\tC}$.

\begin{lem} The weights of $\phi_*\cO_{\tC}$ at the origin are of the form:
$$\mathbb{S}^{(b_1,\dots, b_k)}S_2, \hspace{2cm} 0\leq b_k \leq \dots \leq b_1<k $$
\end{lem}

\begin{proof} $\tC$ is a subvariety in the trivial bundle \eqref{eq.ambienttrivialbundle}. By adjunction, the weights we seek are:
$$\Tor^\bullet\! \big(\cO_{\tC}, \;\cO_{\Fl_1\times \Fl_2}\big) $$
We have described $\tC$ as an iterated generalized complete intersection, so it has a conormal bundle, and we can filter the conormal bundle in the obvious way by line bundles associated to weights of a maximal torus of $\GL(S_1)\times \GL(S_2)$. The result is that our Tor groups can be computed from a spectral sequence that starts with the cohomology groups, on $\Fl_1\times \Fl_2$, of the exterior algebra of the bundle:
 \beq{eq.weightsbundle}\cO^{\oplus k} \oplus \bigoplus_{ i < j } L_i\otimes L_j^{-1} \oplus \bigoplus_{i\neq j} L_i'\otimes L_j^{-1} \oplus \bigoplus_{i < j } L_i'\otimes {L_j'}^{-1}   \eeq
(as above the $L_i$ and $L'_j$ denote the tautological line bundles on $\Fl_1$ and $\Fl_2$). This reduces the problem to Borel-Weil-Bott and some elementary combinatorics. Indeed, given a $\GL(S_2)$-weight
$$(c_1,\dots, c_k)$$
a sufficient condition for the associated line bundle on $\Fl_2$ to have cohomology in our required range is that
$$ c_k, \; c_{k-1} +1, \; \dots, \;  c_1+k-1 \in [0, 2k-2] $$
(the cohomology will be zero if any two of these are equal).  The non-zero $\GL(S_2)$ weights appearing in the bundle \eqref{eq.weightsbundle} are $k-1$ copies of each standard basis vector
$$(0,\dots, 0, 1, 0, \dots, 0)$$
and one each of the negative roots:
$$(0,\dots, 0,-1,0,\dots,1,0,\dots,0)$$

The exterior algebra takes sums of these weights without repetition. Summing the roots gives weights in the range
$$ i-k \leq c_i \leq i-1 $$
and adding the other weights keeps us in the range
$$ i-k \leq c_i \leq k+ i -2$$
as required. 
\end{proof}

\begin{cor}\label{cor.weightsofC} The $\GL(S_2)$ weights of $\cO_C|_0$ are contained in the range:
$$\mathbb{S}^{(b_1,\dots, b_k)}S_2, \hspace{2cm} 0\leq b_k \leq \dots \leq b_1<k $$
\end{cor}

\begin{rem} Part (b) of Proposition \ref{prop.rationality} gives a partial explanation of our weight bounds. The singularity $C$ sits in the ambient vector space \eqref{eq.Z} with codimension $k^2$ and relative equivariant canonical bundle $(\det S_1^\vee \otimes \det S_2)^{k-1}$, hence there is a duality:
$$\Tor^t(\cO_0, \cO_C)^\vee \otimes (\det S_1^\vee \otimes \det S_2)^{k-1} = \Tor^{k^2-t}(\cO_0, \cO_C)$$
This duality swaps our two inequalities $0\leq b_k$ and $b_1<k$. 
\end{rem}

\subsection{Rationality of \texorpdfstring{$C$}{C}}\label{sec.rationality}

Recall that $C$ is an irreducible component of the complete intersection $Z$ \eqref{eq.Z}. 

\begin{lem}\label{lem.jetscheme} The variety $Z$ is semi log canonical and hence $C$ is seminormal. 
\end{lem}

The first statement should be at least plausible because a complete intersection of $k^2$ generic quadrics in $3k^2$ variables has discrepancy $k^2-1$ so is certainly log canonical. The equations for $Z$ are very far from generic, but we can make this heuristic precise by bounding the dimension of the jet scheme of $Z$, following \cite{M}.

\begin{proof} Let $R_m = \C[t]/t^{m+1}$. An $m$-jet in $Z$ is a triple $(p_1^t, a^t, p_2^t)$ of matrices with entries in $R_m$ satisfying $a^t p_2^t = p_1^t a^t$. The jet scheme $J_m(Z)$ is acted on by the algebraic group $\GL_k(R_m)\times \GL_k(R_m)$. Then the proof of Lemma \ref{lem.Zequidimensional} generalizes verbatim to prove that $J_m(Z)$ is equidimensional with:
 $$\dim J_m(Z) = \dim \GL_k(R_m)\times \GL_k(R_m) =  2(m+1)k^2$$
It follows that $Z$ is MJ-log canonical and thus semi log canonical \cite[Prop.~2.13 and Cor.~3.17]{EI}. 

The final implication is \cite[Cor.~1.2]{FL} (\emph{c.f.} \cite[Cor.~6.32]{Kol}) since any Du Bois singularity is seminormal. 
\end{proof}

By choosing an ordering of the eigenvalues of $p_1$ and $p_2$ we can construct a degree $k!$ ramified cover:
$$ B \to C $$
This $B$ is an affine singularity and the cover ramifies over points with repeated eigenvalues. We also have an action of the symmetric group $\mfS_k \acts B$ with $B/\mfS_k= C$. The map $\phi:\tC \to C$ evidently factors via a map $\phi_B: \tC \to B$.

\begin{lem} The fibres of $\phi_B$ are connected.
\end{lem}
\begin{proof} Fix a point in $B$. The fibre of $\phi_B$ over this point is a pair of flags, preserved by $p_1, a, p_2$, with the correct eigenvalues on the quotient lines. Let $E_1\subset S_1$ and $E_2 \subset S_2$ be the eigenspaces of $p_1$ and $p_2$ for the first eigenvalue. To construct a compatible pair of flags the first step is to choose a pair of lines
$$(L_1, L'_1) \in \P E_1\times \P E_2 $$
satisfying $a(L'_1)\subset L_1$; this is a connected variety (though reducible in general). Then the result follows by a simple induction. 
\end{proof}

\begin{lem} $C$ is normal.
\end{lem}

\begin{proof} The variety $\tC$ is normal (it's a generalized complete intersection in a smooth variety with singularities of codimension $>1$) and the map $\phi_B: \tC \to B$ has connected fibres. It follows that the normalization $B^\nu \to B$ of $B$ must be a set-theoretic bijection. 

The $\mfS_k$ action on $B$ must lift to an action on $B^\nu$, and $B^\nu/\mfS_k$ is normal, so it is the normalization $C^\nu$ of $C$. It follows that $C^\nu \to C$ is also a set-theoretic bijection. Since $C$ is seminormal it must therefore be normal.
\end{proof}

\begin{proof}[Proof of Proposition \ref{prop.rationality}.]
(a) Let $Z^\nu$ be the normalization of $Z$ and $\Gamma_Z \subset Z^\nu$ its conductor. Since $Z$ is slc the pair $(Z^\nu, \Gamma_Z)$ is log canonical. 

Away from codimension 2, the intersection of $C$ with the other components of $Z$ is just the component where $\rk a = k-1$, which is (scheme-theoretically) the divisor $\Gamma = (\det a) \subset C$. Since $C$ itself is normal the pair $(C, \Gamma)$ is log canonical. 

The canonical bundle of $Z$ is equivariantly trivial because the determinant of $\Hom(S_2,S_1)$ is the same as the determinant of the ambient vector space \eqref{eq.Z}. So by adjunction we have
\beq{eq.KY}K_C + \Gamma \sim 0 \eeq
and in particular $K_C$ is Cartier. Finally it follows easily from the fact that $(C, \Gamma)$ is log canonical (and $C\setminus \Gamma$ is smooth) that $C$ itself is klt; hence rational. 
\pgap

(b) Rational singularities are Cohen-Macaulay and $K_C$ is Cartier, so $C$ is Gorenstein. From \eqref{eq.KY} we know that $\det a $ is a section of $\omega_C^{-1}$, so $\omega_C= \det S_1^\vee\otimes \det S_2$. 
\end{proof}

\section{Comparison to CKL}\label{sec.comparing}

As before we let  $X_{\pm}$ be the two sides of the stratified Mukai flop, with coordinates $(x,y)$ satisfying the moment map condition $yx=0$ (Section \ref{sec.weights}). Also recall the space $\widetilde{F}$ from Section \ref{sec.Gflop}, which has coordinates:
$$\begin{tikzcd}[] S_1 \ar[yshift=2pt, hook]{r}{x_1}
	&  V
	\ar[yshift=2pt,two heads]{r}{y_2}
	&S_2 \ar[yshift=-2pt, bend left=10]{ll}{a}\end{tikzcd}$$
The fibre product $X_+ \times_{X_0} X_-$ is the subvariety of $\widetilde{F}$ cut out by the equations:
\[ ay_2x_1=0=y_2x_1a\]
It has $k+1$ irreducible components $\mfZ_0,\dots, \mfZ_k$, distinguished by rank conditions:
\[
\mfZ_i = \left\{\, (x_1, y_2, a), \; \rk(a)\leq k-i,\,  \on{null}(y_2 x_1) \geq k-i
\, \right\}
\]
Here we assume that $k\leq n/2$ so that all components are non-empty; if not replace $k$ with $n-k$. 

Let $\widetilde{\mfZ}_i$ be the normalization of $\mfZ_i$, and let $\Theta_i$ be the shifted line bundle:
\[
\Theta_i = \cO_{\widetilde{\mfZ}_i} \otimes (\det S_1)^i \otimes (\det S_2)^{-i} [-i]
\]
Cautis, Kamnitzer and Licata \cite{CKL} construct a \textit{Rickard complex}

\begin{equation}\label{eq:Rickard}
	\Theta = \left\{\, \Theta_k \xlra{d}  \cdots \xlra{d} \Theta_{1} \xlra{d} \Theta_0 \, \right\}
\end{equation}
in $D^b(X_+ \times X_-)$, with $\Theta_0$ in degree $0$. 

\begin{thm}[\cite{CKL}]\label{thm:CKL}
	The complex (\ref{eq:Rickard}) has a unique convolution $\cT = \on{Conv}(\Theta)$, and $\cT$ is the Fourier-Mukai kernel of a derived equivalence
	\[
	\Phi_{\cT} : D^b (T^\vee \Gr(k,V)) \simrightarrow D^b (T^\vee \Gr(V,k))
	\] 
\end{thm}
The convolution $\cT$ is an iterated extension of the sheaves $\Theta_i [i]$. In particular $\cT$ is a sheaf, and it is described explicitly in \cite{Cautis}.

In \cite{Tseu}, the second author transports the complex (\ref{eq:Rickard}) to matrix factorizations.  
Write $$\Psi: D^b (X_{\pm}) \simrightarrow \MF(E_\pm , W)$$
for the Kn\"{o}rrer periodicity equivalences, and set $\Lambda = (\det S_1)^k \otimes (\det S_2)^{-k}$. By \cite[Prop.~3.8]{Tseu}, each functor $\Psi \Phi_{\Theta_i} \Psi^{-1}$ has kernel
\begin{equation*} 
	(\pi_i)_* \cO_{I_i} \otimes \Lambda [-i(i+1)]
\end{equation*}
where $\pi_i : I_i \to E_+ \times E_-$ is defined as follows. 
Let $\mathrm{Fl}(k-i, k, V)$ be the partial flag variety, with tautological bundles $S' \subset S_1\subset V$, and consider the total space
$$\End(S')\oplus \End(S_1) \oplus \End(S_2) \oplus \Hom(S_2, S') \To \mathrm{Fl}(k-i, k, V)\times \Gr(V,k) $$
with coordinates 
\begin{equation*}
	\begin{tikzcd}
		S_1 \arrow[out=120,in=60,loop,looseness=2.7, "p_1"] \arrow[r, hook, "x_1"] 
		& V \arrow[r,  two heads, "y_2"] 
		& S_2 \arrow[out=120,in=60,loop,looseness=2.7, "p_2"] \arrow[ld,"a_2"]  \\
		& S' \arrow[lu,hook,"a_1"] \arrow[out=300,in=240,loop,looseness=2.7, "p'"]  &
	\end{tikzcd}
\end{equation*}
where  $x_1,y_2, a_1$ have full rank. Then $I_i$ is the subvariety 
\begin{equation*}
	I_i = \{ a_1 p' = p_1 a_1, \;\; a_2 p_2 = p' a_2 \}
\end{equation*}
and $\pi_i$ sends a point to $(x_1, y_2, a_1 a_2, p_1, p_2)$. Since $a_1$ is injective, the first equation says that $S'$ is $p_1$-invariant and $p' = p_1 |_{S'}$. 
The differentials of (\ref{eq:Rickard}) are transported to morphisms of matrix factorizations, still denoted by $d$. 
Omitting the common factor $\Lambda$, we obtain a complex in $\MF(E_+ \times E_-, W' - W)$: 
\begin{equation}\label{eq:tRic}
	(\pi_k)_* \cO_{I_k} [-k^2 -k] \; \xlra{d}\; \cdots\; \xlra{d} \;(\pi_1)_* \cO_{I_1} [-2]\; \xlra{d} \; \cO_{I_0}
\end{equation}
Its convolution is unique, because the convolution of (\ref{eq:Rickard}) is unique and $\Psi$ is an equivalence. After tensoring with $\Lambda$ it is the kernel of $\Psi \Phi_{\cT} \Psi^{-1}$. 

\begin{thm}[{\cite[Thm.~3.3]{Tseu}}]
 The convolution of \eqref{eq:tRic} is the kernel for the window equivalence corresponding to the window
	$$\big\{ \mathbb{S}^{(a_1,\dots, a_k)} S, \; 0\leq a_k \leq \cdots \leq a_1 <n \, \big\}$$ 
\end{thm}

Note that this is exactly the window $\mfW'$ (\ref{eq:mfW'}). Hence by Theorem \ref{thm.main}:

\begin{cor}\label{cor.convolution}
The convolution of (\ref{eq:tRic}) is isomorphic to $\cO_{\overline{\Delta}}$ in: $$\MF(E_+\times E_-, W'-W)$$
\end{cor} 

Our goal for the remainder of this section is to unpack Corollary~\ref{cor.convolution} heuristically.
We will argue that (\ref{eq:tRic}) should be a resolution of $\cO_{\overline{\Delta}}$ which, restricted to each rank stratum of $I_0$, becomes a Koszul complex. 

\pgap
Recall the variety $D$ from (\ref{eq.D}), which is a vector bundle over $\Gr(k,V) \times \Gr(V,k)$ with fibre $\Hom(S_2,S_1)\oplus \End(S_1) \oplus \End(S_2)$. We have $I_0 = \{p_1 a = ap_2\} \subset D$, so $I_0$ is a fibre bundle over $\Gr(k,V) \times \Gr(V,k)$ with fibre the variety $Z$ of Section \ref{sec.weightsofC}. 
It has $k+1$ irreducible components
\[
C_i = \overline{\{\rk(a) = k-i \} } \subset I_0, \quad i=0,\dots, k
\]
and $C_0 = \overline{\Delta}$ by Section \ref{sec.geometry}.  
We stratify $I_0$ by the rank of $a$: let
\[
I_0^{(j)} = \{\rk(a) = k-j\} \subset I_0, \quad j=0,\dots, k
\] 
On $I_0^{(j)}$, the map $a$ has kernel $K\subset S_2$ and image $A\subset S_1$, of ranks $j$ and $k-j$. The equation $p_1 a = a p_2$ says that $p_2$ preserves $K$, $p_1$ preserves $A$, and $a$ intertwines them on $S_2 /K \cong A$. So the characteristic polynomials of $p_1$ and $p_2$ share the factor of degree $k-j$ coming from $A$. Write 
\[
q_1 \in \End(S_1/A), \quad q_2 = p_2|_{K} \in \End(K)
\]
for the induced endomorphisms, then 
\[
\overline{\Delta} \cap I_0^{(j)} = \{\chi_{q_1} = \chi_{q_2}\}
\]
This locus is cut out by the $j$ coefficients of $\chi_{q_1} - \chi_{q_2}$, which have R-charges $2,4,\dots,2j$ and form a regular sequence. So, $\cO_{\overline{\Delta}}|_{I_0^{(j)}}$ is resolved by the Koszul complex
$\Wedge^\bullet F_j$, where:
\vspace{-14pt}
$$F_j = \bigoplus_{r=1}^j \cO [-2r]$$
On the other hand, over $I_0^{(j)}$ the condition $a(S_2) \subset S'$ becomes $A\subset S'$. So $\pi_i^{-1}(I_0^{(j)})$ is the partial Grothendieck-Springer resolution of $q_1$-invariant subspaces of dimension $j-i$ in $S_1/A$, and is empty for $i>j$. 
Its pushforward to $I_0^{(j)}$ is free, with fibre $H^\bullet (\Gr(j-i,j))$ over nilpotent $q_1$ and graded rank the Gaussian binomial $\binom{j}{i}_q$, where $q$ counts R-charge $2$ (\emph{c.f.} \cite[Eq.~(26)]{Tseu}). The graded rank of $\Wedge^i F_j$ is
$$e_i(q,q^2,\dots,q^j) = q^{i(i+1)/2}\binom{j}{i}_q $$
So on every stratum the $i$-th term of \eqref{eq:tRic}, including its shift $[-i(i+1)]$, has exactly the graded rank of the $i$-th term of the Koszul complex. The resulting picture is Figure \ref{fig.stratifiedKoszul}.

\begin{figure}[ht]
	\centering
	\small
	\begin{tikzcd}[row sep=1.4em, column sep=1.1em]
		& I_0^{(0)} & I_0^{(1)} & I_0^{(2)} & \cdots & I_0^{(k)} \\[-0.6em]
		(\pi_k)_*\cO_{I_k}[-k^2-k] \ar[d] & 0 & 0 & 0 & & \Wedge^k F_k \ar[d] \\
		\vdots \ar[d] & \vdots & \vdots & \vdots & & \vdots \ar[d] \\
		(\pi_2)_*\cO_{I_2}[-6] \ar[d,"d"] & 0 & 0 & \Wedge^2 F_2 \ar[d] & & \Wedge^2 F_k \ar[d] \\
		(\pi_1)_*\cO_{I_1}[-2] \ar[d,"d"] & 0 & F_1 \ar[d] & F_2 \ar[d] & & F_k \ar[d] \\
		\cO_{I_0} \ar[d] & \cO \ar[d] & \cO \ar[d] & \cO \ar[d] & & \cO \ar[d] \\
		\cO_{\overline{\Delta}} & \cO & \cO_{\overline{\Delta}}|_{I_0^{(1)}} & \cO_{\overline{\Delta}}|_{I_0^{(2)}} & & \cO_{\overline{\Delta}}|_{I_0^{(k)}}
	\end{tikzcd}
	\caption{The complex \eqref{eq:tRic} as a `stratified' Koszul complex: the entry in row $i$ and column $j$ is the (underived) restriction of the $i$-th term to the stratum $I_0^{(j)}$, and each column is a Koszul resolution.}
	\label{fig.stratifiedKoszul}
\end{figure}

\bibliographystyle{halphanum}

\end{document}